\documentclass[11pt]{article}
\usepackage{epsfig }
\usepackage{latexsym}
\usepackage{amsfonts,amsmath,amssymb,textcomp}
\allowdisplaybreaks

\newtheorem{thm}{Theorem}[section]
\newtheorem{lem}[thm]{Lemma}
\newtheorem{cor}[thm]{Corollary}
\newtheorem{prop}[thm]{Proposition}
\newtheorem{rem}[thm]{Remark}
\newtheorem{deff}[thm]{Definition}
\newtheorem{conj}[thm]{Conjecture}
\newtheorem{key}[thm]{Keywords}
\newtheorem{prob}[thm]{Problem}
\newenvironment{probb}{ \begin{prob} \rm}{ \end{prob} }
\newcommand{\bprobb}{\begin{probb}}
\newcommand{\eprobb}{\end{probb}}
\newcommand{\bth}{\begin{thm}}
\newcommand{\ethGL}{\end{thm}}
\newcommand{\bconj}{\begin{conj}}
\newcommand{\econj}{\end{conj}}
\newcommand{\bkey}{\begin{key}}
\newcommand{\ekey}{\end{key}}
\newcommand{\bl}{\begin{lem}}
\newcommand{\el}{\end{lem}}
\newcommand{\bdeff}{\begin{deff}}
\newcommand{\edeff}{\end{deff}}
\newcommand{\bcor}{\begin{cor}}
\newcommand{\ecor}{\end{cor}}
\newcommand{\bprop}{\begin{prop}}
\newcommand{\eprop}{\end{prop}}
\newcommand{\brem}{\begin{rem}}
\newcommand{\erem}{\end{rem}}
\newcommand{\beq}{\begin{equation}}
\newcommand{\eeq}{\end{equation}}
\newcommand{\beqn}{\begin{eqnarray}}
\newcommand{\eeqn}{\end{eqnarray}}
\newcommand{\beqns}{\begin{eqnarray*}}
\newcommand{\eeqns}{\end{eqnarray*}}
\newcommand{\ba}{\begin{array}}
\newcommand{\ea}{\end{array}}
\newcommand{\bit}{\begin{itemize}}
\newcommand{\eit}{\end{itemize}}
\newcommand{\ben}{\begin{enumerate}}
\newcommand{\een}{\end{enumerate}}

\newcommand{\BO}{\mathcal{O}}

\newcommand{\non}{\nonumber}
\newcommand{\babs}{\begin{abstract}}
\newcommand{\eabs}{\end{abstract}}
\newcommand{\bal}{\begin{align}}
\newcommand{\bals}{\begin{align*}}
\newcommand{\bs}{\begin{skip}}
\newcommand{\eal}{\end{align}}
\newcommand{\eals}{\end{align*}}
\newcommand{\es}{\end{skip}}

\newcommand{\ra}{\rightarrow}

\newcommand{\lp}{\left (}
\newcommand{\rp}{\right )}
\newcommand{\lb}{\left [}
\newcommand{\rb}{\right ]}

\newcommand{\al}{\alpha}

\newcommand{\eps}{\varepsilon}
\newcommand{\II}{\infty}

\newcommand{\gam}{\gamma}
\newcommand{\Gam}{\Gamma}

\def\1{{\ifmmode 1\mskip-1.5\thinmuskip\mathrm{l}%
        \else\textrm{1\hskip -.23em l}\fi}}

\newcommand{\bin}[2]
{
{#1\choose #2}
}

\title{\bf Two applications of  polylog functions and Euler sums }
\author{Guy Louchard\thanks{Universit\'e Libre de Bruxelles,
D\'epartement d'Informatique, CP 212, Boulevard du Triomphe, B-1050
Bruxelles, Belgium, email: louchard@ulb.ac.be}
}
\date{\today}
\begin{document}
\maketitle

\babs
 Let
$I(n):=\int_0^1 [x^n+(1-x)^n]^\frac1n dx.$ In this paper, we show that $I(n)= \sum_0^\II \frac{I_i}{n^i},n\ra\II$ and we compute
$I_i, i =0..5$, obtained by polylog functions and Euler sums. As a corollary, we obtain explicit expressions for some integrals involving functions $ u^i, exp(-u),
 (1 +exp(-u))^j , ln(1 + exp(-u))^k$ . As another asymptotic result, let  $S_0(z):=\frac{Li_m(1)}{Li_m(1)-Li_m(z)}$, where $Li_m(z)$ is the polylog function.
 We provide the asymptotic behaviour of $S_n,n\ra\II$ where $S_n:=[z^n]S_0(z)$.
This paper fits within the framework of analytic combinatorics.
\eabs


\section{Introduction}\label{S1}
Some time ago,  the following question  was circulating among the Mathematical problems aficionados: let

\[I(n):=\int_0^1 [x^n+(1-x)^n]^\frac1n dx,\]
what are  $I_0:=\lim I(n), I_2:= \lim n^2(I(n)-I_0), n\ra\II$ ? I found it interesting to look at a deeper asymptotic analysis 
of $I(n)$ and found actually that, asymptotically, 
\[I(n)= \sum_0^\II \frac{I_i}{n^i},n\ra\II,\]
where $I_i$ are curiously obtained by polylog functions and Euler sums. In this paper we compute $I_i, i =
0..5$. As a corollary, we obtain explicit expressions for some integrals involving functions $ u^i, exp(-u),
 (1 +exp(-u))^j , ln(1 + exp(-u))^k$. About polylog functions, see de Doelder, \cite{DO91},
Apostol, \cite{AP10}, Lewin, \cite{LE81}, and about
Euler sums, see Flajolet, Salvy, \cite{FlSa98}, Xu, \cite{XU17}.

Another problem arose in  some work in progress on dynamical systems by G\'{o}mez-Aiza and Ward \cite{GOWA17}. Ward asked the following question: the polylog function is defined as

\[Li_m(z):=\sum_1^\II \frac{z^n}{n^m}.\]
Set
\[ S_0(z):=\frac{Li_m(1)}{Li_m(1)-Li_m(z)}\]
and
\[S_n:=[z^n]S_0(z) .\]
What is the asymptotic behaviour of $S_n,n\ra\II$? In this paper, we provide the asymptotics of $S_n,m=3,4$, up to the $1/n^3$
term. Next terms can be mechanically computed.

\section{A first analysis of $I_n$}\label{S2}
We have

\[I(n)=2\int_0^{1/2} [x^n+(1-x)^n]^\frac1n dx=2\int_0^{1/2} (1-x)\lb 1+\lp\frac{x}{1-x}\rp^n\rb ^\frac1n dx,\]
and

 \bals
&0\leq\frac{x}{1-x}\leq 1,\mbox{ let }\\
F(n)&:=\lb 1+\lp\frac{x}{1-x}\rp^n\rb ^\frac1n 
\sim \exp\lb \left.  \lp\frac{x}{1-x}\rp^n\right/ n \rb,\\
&\left. \lp\frac{x}{1-x}\rp^n\right/n \ra 0,n\ra\II,\mbox{ exponentially if } x<1/2, \mbox{ as } 1/n, \mbox{ if } x=1/2.
\end{align*}
Hence the asymptotic behaviour of $I(n)$ is related to the behaviour of $F(n)$ in the neighbourhood of $x=1/2$. We set $x=1/2-y$
and get 

\[I_0=2\int_0^{1/2} (1-x)dx=\frac34.\]
We now expand $I_n$ up to the $1/n^5$ term.
\bals
I(n)&= 2\int_0^{1/2} \lp\frac12+y\rp\lb 1+\lp \frac{1-2y}{1+2y} \rp^n\rb ^\frac1n dy\\
&= 2\int_0^{1/2} \lp\frac12+y\rp\lb 1+\lp  1-4y+8y^2-16y^3+32y^4+\BO(y^5) \rp^n\rb ^\frac1n dy,\mbox{ and with }y=\frac{u}{4n},\\
&= \frac{2}{4n}\int_0^{2n} \lp\frac12+\frac{u}{4n}\rp\lb 1+\lp   1-u/n+1/2u^2/n^2-1/4u^3/n^3+1/8u^4/n^4+\BO(u^5/n^5)\rp^n\rb ^\frac1n du\\
&=\frac{2}{4n}\int_0^{2n} \lp\frac12+\frac{u}{4n}\rp\lb 1+ \exp\lp -u-1/12 u^3/n^2 +\BO(u^5/n^4) \rp  \rb^\frac1n du\\
&= \frac{2}{4n}\int_0^{2n} \lp\frac12+\frac{u}{4n}\rp
 \lb 1+ \exp(-u)-1/12\exp(-u)u^3/n^2+\exp(-u)\BO(u^5/n^4) \rb ^\frac1n du\\
&= \frac{2}{4n}\int_0^{2n} \lp\frac12+\frac{u}{4n}\rp
 \exp\lb \ln(1+\exp(-u))/n-1/12\exp(-u)u^3/((1+\exp(-u))n^3)+\exp(-u)\BO(u^5/n^5)\rb du\\
&= \frac{2}{4n}\int_0^{2n}  \lp\frac12+\frac{u}{4n}\rp
  \bigg[ 1+\ln(1+\exp(-u))/n+1/2\ln(1+\exp(-u))^2/n^2  \\
&+1/12\Bigl[ -\exp(-u)u^3/(1+\exp(-u))+2\ln(1+\exp(-u))^3\Bigr]/n^3\\
& +1/24\ln(1+\exp(-u))\Bigl[-2\exp(-u)u^3/(1+\exp(-u))+\ln(1+\exp(-u))^3\Bigr]/n^4\\
&+\exp(-u)\BO(u^5/n^5)\bigg] du
\end{align*}
\newpage
\bals
&=\frac34+ \int_0^\II\lb \rule{0mm}{7mm} 1/4\ln(1+\exp(-u))/n^2+\Bigl[ 1/8u\ln(1+\exp(-u))+1/8\ln(1+\exp(-u))^2\Bigr]/n^3\right.\\
&+\bigg[ 1/16u\ln(1+\exp(-u))^2
 +1/48\Bigl[-\exp(-u)u^3/(1+\exp(-u))+2\ln(1+\exp(-u))^3\Bigr]\bigg]/n^4\\
& +1/96\bigg[ -\exp(-u)u^4/(1+\exp(-u))+2u\ln(1+\exp(-u))^3\\
&  -2\ln(1+\exp(-u))\exp(-u)u^3/(1+\exp(-u))+\ln(1+\exp(-u))^4\bigg] /n^5\\
&\left.+\exp(-u)\BO(u^6/n^6)\rule{0mm}{7mm} \rb du.
\end{align*}
We immediately recover $I_0$. The computation of $I_i,i\geq 1$ is detailed in the next sections.
\section{Computation of $I_1,I_2,I_3$}\label{S3}
We have
\[I_1=0,\]
\[I_2= \int_0^\II  1/4\ln(1+\exp(-u))du=1/4\sum_1^\II \frac{(-1)^{i+1}}{i^2}=\frac{\pi^2}{48},\]
\[I_3=\int_0^\II [1/8 u\ln(1+\exp(-u))+1/8\ln(1+\exp(-u))^2]du=11/32\zeta(3)+1/8 I\ln(2)^2\pi+1/4\ln(2)Li_2(2)-1/4Li_3(2),\]
where the polylog function is defined by
\[Li_n(z)=\sum_1^\II \frac{z^k}{k^n}.\]
But we know that
\[Li_n(z)=-(-1)^n Li_n(1/z)-\frac{(2\pi I)^n}{n!}B_n\lp\frac12+\frac{\ln(-z)}{2\pi I} \rp,z\notin [0,1],\] 
where $B_n(x)$ is the $n$th Bernoulli polynomial,
and
\bals
Li_2(1/2)&=1/12\pi^2-1/2\ln(2)^2,  \mbox{  hence }Li_2(2)= 1/4\pi^2-I\pi\ln(2),\\
Li_3(1/2)&=7/8\zeta(3)-1/12\pi^2\ln(2)+1/6\ln(2)^3,  \mbox{  hence }Li_3(2)=7/8\zeta(3)+1/4\pi^2\ln(2)-1/2 I\pi\ln(2)^2.
\end{align*}
The values of $Li_k(1/2),k\geq 4$ are not known to be related to classical constants.

This leads to

\[I_3=\frac{\zeta(3)}{8}.\]
Another, more elegant, way to compute $I_3$ is to turn to Euler sums. Following  Flajolet, Salvy, \cite{FlSa98}, we have
\bals
S_{p,q}^{+-}&:=\sum_{k=1}^\II (-1)^{k-1}\frac{H_k^{(p)}}{k^q},H_n^{(p)}:=\sum_{j=1}^n \frac{1}{j^p},\\
\overline{\zeta}(s)&:=(1-2^{1-s})\zeta(s),\overline{\zeta}(1):=\ln(2),\\
2S_{1,q}^{+-}&=(q+1)\overline{\zeta}(q+1)-\zeta(q+1)-2\sum_{k=1}^{q/2-1}\overline{\zeta}(k)\zeta(q+1-2k),1+q\mbox{ odd },\\
2S_{1,2}^{+-}&=\frac{5}{4}\zeta(3),\\
2S_{1,4}^{+-}&=59/16\zeta(5)-1/6\pi^2\zeta(3).
\end{align*}
Hence
\bals
\int_0^\II  u\ln(1+\exp(-u))du &=\sum_1^\II \frac{(-1)^{i+1}}{i^3}=\frac34\zeta(3),\\
\int_0^\II \ln(1+\exp(-u))^2 du &= \sum_1^\II \sum_1^\II \frac{(-1)^{i+j}}{ij(i+j)}
=\sum_{k=2}^\II\sum_{i=1}^{k-1}\frac{(-1)^k}{ki(k-i)}
=\sum_{k=2}^\II \frac{(-1)^k}{k^2}2H_{k-1}\\
&=2\sum_{k=1}^\II \frac{(-1)^k}{k^2}H_{k}-2\sum_{k=1}^\II \frac{(-1)^k}{k^3}
=-2S_{1,2}^{+-}+\frac32\zeta(3)=\frac{\zeta(3)}{4}.\\
\end{align*}
We immediately recover $I_3$.

Another similar sum will be used in the next section: for $p+q$ odd,
\bals
S_{p,q}^{+-}&:=\lb  (1-(-1)^p)\zeta(p)\overline{\zeta}(q)+\overline{\zeta}(p+q)
+2\sum_{k=0}^{\lfloor p/2 \rfloor}\bin{q+p-2k-1}{q-1}(-1)^{p-2k+1}\overline{\zeta}(q+p-2k)\overline{\zeta}(2k) \right. \\
&\left. +2(-1)^p \sum_{0}^{\lfloor q/2 \rfloor}\bin{p+q-2k-1}{p-1}\zeta(p+q-2k)\overline{\zeta}(2k)\rb\left/\rule{0mm}{20mm}2 \right.,\\
S_{2,3}^{+-}&= -11/32\zeta(5)+5/48\zeta(3)\pi^2.\\
\end{align*}

\section{Computation of $I_4$}\label{S4}
Now we have
\bals
S_1&:=\int_0^\II\frac{u^3 e^{-u}}{1+e^{-u}}du =\frac{7\pi^4}{120},\\
S_2&:=\int_0^\II \ln(1+\exp(-u))^3 du  =\ln(2)^3\pi I+3\ln(2)^2Li_2(2)-6\ln(2)Li_3(2)+6Li_4(2)-1/15\pi^4\\
&=1/4\pi^2\ln(2)^2-21/4\ln(2)\zeta(3)-6Li_4(1/2)+1/15\pi^4-1/4\ln(2)^4,\\
S_3&:= \int_0^\II u\ln(1+\exp(-u))^2 du=\sum_{k=2}^\II \frac{(-1)^k}{k^3}2H_{k-1}
=2\sum_{k=1}^\II \frac{(-1)^k}{k^3}H_{k}-2\sum_{k=1}^\II \frac{(-1)^k}{k^4}
=-2S^{+-}_{1,3}-2\lp  -\frac{7}{720}\pi^4\rp,\\
S^{+-}_{1,3}&=-2Li_4(1/2)+11/4\zeta(4)+1/2\zeta(2)\ln(2)^2-1/12\ln(2)^4-7/4\zeta(3)\ln(2),
\mbox{ this is }\mu_1 \mbox{ in } \cite{FlSa98},\\
S_3&= 4Li_4(1/2)-1/24\pi^4-1/6\pi^2\ln(2)^2+1/6\ln(2)^4+7/2\ln(2)\zeta(3).
\end{align*}
Hence
\[I_4:=\frac{1}{48}\lb -S_1+2S_2+3S_3\rb=-\frac{\pi^4}{960}.\]
\section{Computation of $I_5$}\label{S5}
We compute

\[S_4:=\int_0^\II \frac{u^4 e^{-u}}{1+e^{-u}}du=45/2\zeta(5).\]
Now we turn to $S_5$, which is the most intricate case of our integral expressions: 
\bals
S_5&:=\int_0^\II u\ln(1+\exp(-u))^3 du=\sum_{i=1}^\II \sum_{j=1}^\II \sum_{\ell=1}^\II\frac{(-1)^{i+j+\ell+1}}{ij\ell(i+j+\ell)^2}=
\sum_{k=3}^\II\frac{(-1)^{k+1}}{k^2}\sum_{v=2}^{k-1}\frac{1}{k-v}\sum_{i=1}^{v-1}\frac{1}{i(v-i)}\\
&=\sum_{v=2}^{\II}2\frac{H_{v-1}}{v}\sum_{k=v+1}^\II\frac{(-1)^{k+1}}{k^2(k-v)}\\
&=2\sum_{v=2}^{\II}\frac{(-1)^{k-1}}{k^2}\sum_{j=1}^{k-1}\frac{H_{j-1}}{j(k-j)}\\
&=2\sum_{v=2}^{\II}\frac{(-1)^{k-1}}{k^3}\sum_{j=1}^{k-1}H_{j-1}\Big[\frac1j+\frac1{k-j}\Big]\\
&=\sum_{v=2}^{\II}\frac{(-1)^{k-1}}{k^3}\Big[H_{k-1}^2-H^{(2)}_{k-1}+2H_{k-1}^2-2H^{(2)}_{k-1}\Big]\\
&=3\sum_{v=2}^{\II}\frac{(-1)^{k-1}}{k^3}\Big[H_{k-1}^2-H^{(2)}_{k-1}\Big].
\end{align*}
But
\bals
H_{n-1}^2/n^3&=[H_n^2-2H_n/n+1/n^2]/n^3,\\
H_{n-1}^{(2)}/n^3&=[H_{n}^{(2)}-1/n^2]/n^3.
\end{align*}
Hence
\bals
S_5&=-3(T_1+T_2), \mbox{ with }\\
T_2&=-\lb -S^{+-}_{2,3} +15/16\zeta(5)\rb=
                5/48 \zeta(3) \pi^2-\frac{41}{32}  \zeta(5),\\
T_1&=T_3+2 S^{+-}_{1,4}-15/16\zeta(5),\\
T_3&=\sum_{k=1}^\II (-1)^{k}\frac{H_k^{2}}{k^3}\\
&=-(4Li_5(1/2)+4\ln(2)Li_4(1/2)+2/15\ln(2)^5+7/4\zeta(3)\ln(2)^2\\
&-19/32\zeta(5)-2/3\zeta(2)\ln(2)^3 -11/8\zeta(2)\zeta(3)),
\mbox{ see \cite{XU17}, where many recent references can be found, so}\\
T_1&= -4 Li_5( 1/2) - 4 \ln(2) Li_4( 1/2) - 2/15 \ln(2)^5
          - 7/4 \zeta(3) \ln(2)^2+\frac{107}{32}   \zeta(5) + 1/9 \pi^2  \ln(2)^3
         + 1/16 \zeta(3) \pi^2, \\
				\mbox{ and finally} &\\
S_5&= 12 Li_5 1/2) + 12 \ln(2) Li_4( 1/2) + 2/5 \ln(2)^5
         + 21/4 \zeta(3) \ln(2)^2-\frac{99}{16}  \zeta(5) - 1/3 \pi^2  ln(2)^3
         - 1/2 \zeta(3) \pi^2,\\                                           
S_6&:=\int_0^\II \frac{u^3 e^{-u} \ln(1+\exp(-u))}{1+e^{-u}}du=\int_0^\II u^3\sum_{k=2}^\II e^{-uk}(-1)^k\sum_{i=1}^{k-1}\frac1i du=
3!\sum_{k=2}^\II \frac{1}{k^4}(-1)^k H_{k-1}\\
&=3!\lb \sum_{k=1}^\II \frac{1}{k^4}(-1)^k H_{k}-\sum_{k=1}^\II \frac{1}{k^5}(-1)^k\rb =
3!\lb -S^{+-}_{1,4}+\frac{15}{16}\zeta(5)\rb =-87/16\zeta(5)+1/2\pi^2\zeta(3),\\
S_7&:=\int_0^\II\ln(1+\exp(-u))^4 du \\
&= 2/3\pi^2\ln(2)^3-21/2\ln(2)^2\zeta(3)-24\ln(2)Li_4(1/2)-4/5\ln(2)^5-24Li_5(1/2)+24\zeta(5).
 \end{align*}
Hence
\[I_5=\frac{1}{96}
\lb -S_4+2S_5-2S_6+S_7\rb=-1/48\zeta(3)\pi^2.\]
Let us summarize our results in the following theorem:
\bth
Let
\[I(n):=\int_0^1 [x^n+(1-x)^n]^\frac1n dx.\]
We have
\[I(n)= \sum_0^\II \frac{I_i}{n^i},n\ra\II,I_0=\frac34,I_1=0,I_2=\frac{\pi^2}{48},I_3=\frac{\zeta(3)}{8},I_4=-\frac{\pi^4}{960},I_5=-1/48\zeta(3)\pi^2.\]
\ethGL
We also have the following corollary (many other similar results can be found in Xu \cite{XU17})
\bcor
\bals
&\int_0^\II  1/4\ln(1+\exp(-u))du=\frac{\pi^2}{48},\\
&\int_0^\II  u\ln(1+\exp(-u))du =\frac34\zeta(3),\\
&\int_0^\II \ln(1+\exp(-u))^2 du = \frac{\zeta(3)}{4},\\
&\int_0^\II\frac{u^3 e^{-u}}{1+e^{-u}}du =\frac{7\pi^4}{120},\\
&\int_0^\II \ln(1+\exp(-u))^3 du 
=1/4\pi^2\ln(2)^2-21/4\ln(2)\zeta(3)-6Li_4(1/2)+1/15\pi^4-1/4\ln(2)^4,\\
&\int_0^\II u\ln(1+\exp(-u))^2 du=
 4Li_4(1/2)-1/24\pi^4-1/6\pi^2\ln(2)^2+1/6\ln(2)^4+7/2\ln(2)\zeta(3),\\
&\int_0^\II \frac{u^4 e^{-u}}{1+e^{-u}}du=45/2\zeta(5),\\
&\int_0^\II u\ln(1+\exp(-u))^3 du\\
&= 12 Li_5 1/2) + 12 \ln(2) Li_4( 1/2) + 2/5 \ln(2)^5
         + 21/4 \zeta(3) \ln(2)^2-\frac{99}{16}  \zeta(5) - 1/3 \pi^2  ln(2)^3
         - 1/2 \zeta(3) \pi^2,\\                                           
&\int_0^\II \frac{u^3 e^{-u} \ln(1+\exp(-u))}{1+e^{-u}}du=
-87/16\zeta(5)+1/2\pi^2\zeta(3),\\
&\int_0^\II\ln(1+\exp(-u))^4 du \\
&= 2/3\pi^2\ln(2)^3-21/2\ln(2)^2\zeta(3)-24\ln(2)Li_4(1/2)-4/5\ln(2)^5-24Li_5(1/2)+24\zeta(5).
\end{align*}
\ecor
We have two Open problem:

Open problem $1$: how to explain the relatively simple $I_i$ expressions?

Open problem $2$: can we find `easily' similar computations for $I_i,i\geq 6$?

\begin{figure}[htbp]
	\centering
		\includegraphics[width=0.8\textwidth,angle=0]{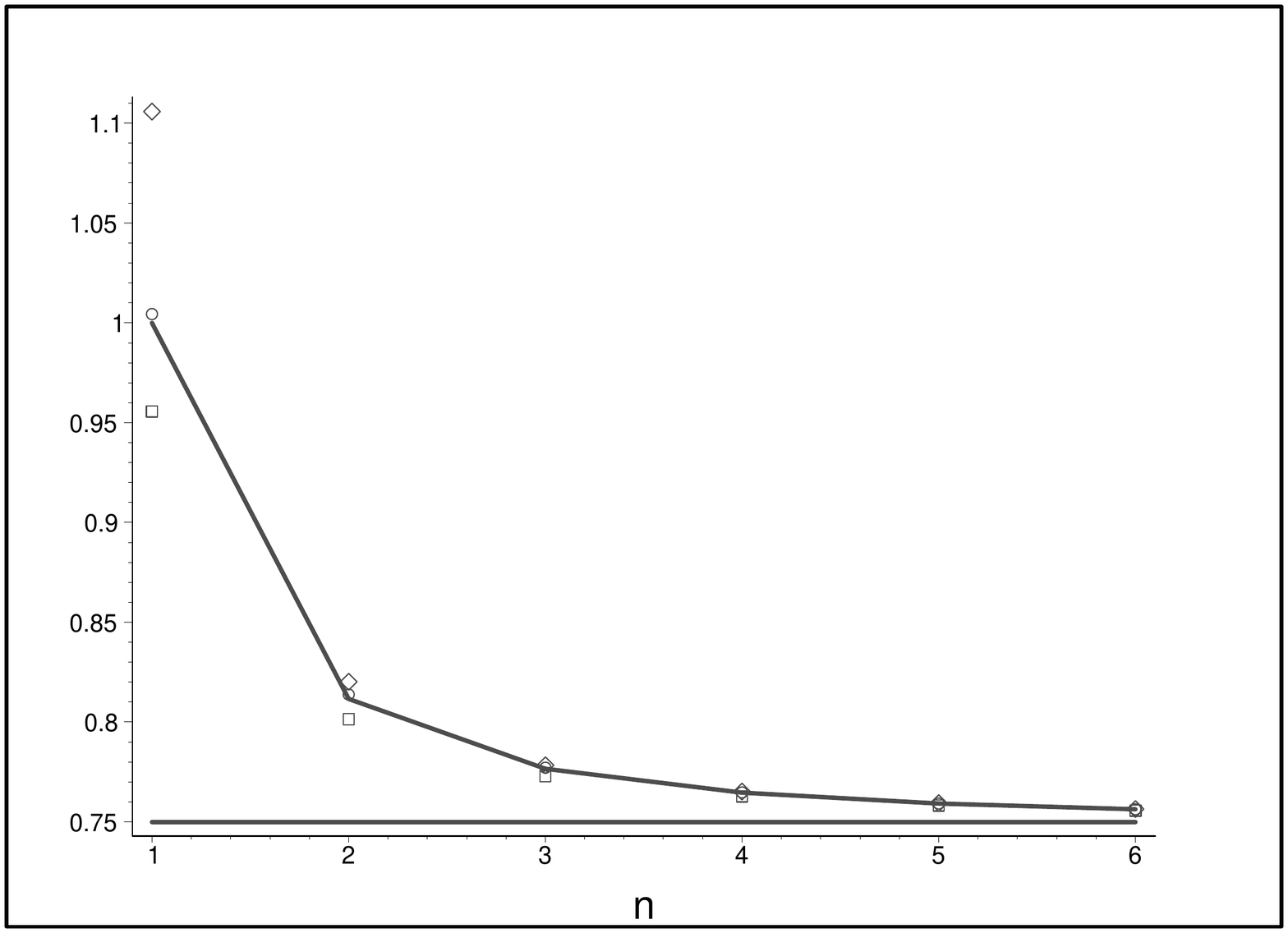}
	\caption{ $I(n)$(line),$I_0$(line),$I_{n,2}$(box),$I_{n,3}$(diamond),$I_{n,4}$(circle)}
	\label{F1}
\end{figure}

Let
\bals
I_{n,2}&:=I_0+I_2/n^2,\\
I_{n,3}&:=I_0+I_2/n^2+I_3/n^3,\\
I_{n,4}&:=I_0+I_2/n^2+I_3/n^3+I_4/n^4.\\
\end{align*}
\newpage
To check the quality of our asymptotics, we display, in Figure \ref{F1}, $I_0,I(n),I_{n,2},I_{n,3},I_{n,4}$.

\section{A first analysis of $S_n$}\label{S6}
\bal 
&\mbox{We know that}           \non\\
 Li_m(z)&=\frac{(-1)^m}{(m-1)!}w^{m-1}(\ln(w)-H_{m-1})+\sum_{j\neq m-1,j\geq 0}\frac{(-1)^j}{j!}\zeta(m-j)w^j,\mbox{ see \cite{FlSe09},VI.20, with }                                 \label{E1}\\
 w&:=-\ln(z).           \non
\end{align}
The singularity $z=1$  in  $S_0(z)$  leads to the desired expansion of $S_n$.  Set
\bals
   L&:=\ln\lp \frac{1}{1-z} \rp=Li_1(z),\\
 L_{k,n}&:=[z^n]L^k,\mbox{ with}\\
L_{1,n}&=\frac1n.\\
\mbox{Let}&\\
\eps&:=1-z,\\
S(\eps,L)&:=\mbox{ expansion of }  S_0(1-\eps)      \mbox{ w.r.t. }\eps,\\
&\mbox{ we compute successively}\\
 D_{i,j}&:=[\eps^i L^j]S(\eps,L),\mbox{ which depends on }m,\\
G_{i,j,n}&:=[z^n]\eps^i L^j,\mbox{ which is independent of }m,\\
T_{i,j,n}&= D_{i,j}G_{i,j,n}=[z^n\eps^i L^j]S(\eps,L),\\
S_n&:=[z^n]S_0(z) ,\\
T_n&:=\mbox{ asymptotics of }S_n =\sum_i\sum_j T_{i,j,n},n\ra\II.\\
&\mbox{ We also define}\\
C_{n,k}&:=\mbox{ asymptotics of }S_n \mbox{ up to the }1/n^k \mbox{ term}.
\end{align*}
In this paper, we will compute $C_{n,k}, k=0..3,m=3,4$, but as we will see, more terms can be mechanically obtained. We will also show some graphs of $C_{n,k}$.

\section{Some asymptotics for $L_{k,n}$}\label{S7}
\bals
&\mbox{ Set}\\
H &:= \Gam(n+\al)/(\Gam(\al)\Gam(n+1)),\mbox{ see [1],VI.7},\\
\frac{\partial^2H(\al)}{\partial \al^2}&=[\psi(1,n+\al)+\psi(n+\al)^2
-2\psi(n+\al)\psi(\al)
+\psi(\al)^2-\psi(1,\al)]\Gam(n+\al)/(\Gam(\al)\Gam(n+1)),\\
&\mbox{ where }  \psi(n,x) \mbox{ is the nth polygamma function, which is the nth derivative of the digamma function},\\
&\mbox{ we have}\\
L_{2,n}&=\lim_{\al\ra 0}\frac{\partial^2H(\al)}{\partial \al^2}=(2\psi(n)+2\gam)/n,\\
L_{2,n}&=(2\ln(n)+2\gam)/n-1/n^2-1/(6n^3)+1/(60n^5)+\BO(1/n^6),\mbox{ see [1], Figure VI.5 for the first terms},\\
\frac{\partial^3H(\al)}{\partial \al^3}&=[\psi(2,n+\al)+3\psi(1,n+\al)\psi(n+\al)-3\psi(1,n+\al)\psi(\al)\\
&+\psi(n+\al)^3-3\psi(n+\al)^2\psi(\al)+3\psi(n+\al)\psi(\al)^2\\
&-3\psi(n+\al)\psi(1,\al)-\psi(\al)^3+3\psi(\al)\psi(1,\al)-\psi(2,\al)]\Gam(n+\al)/(\Gam(\al)\Gam(n+1)),\\
L_{3,n}&=\lim_{\al\ra 0}\frac{\partial^3H(\al)}{\partial \al^3}=1/2(6\psi(1,n)+12\psi(n)\gam-\pi^2+6\gam^2+6\psi(n)^2)/n,\\
L_{3,n}&=(3\ln(n)^2+6\ln(n)\gam-1/2\pi^2+3\gam^2)/n+(3-3\ln(n)-3\gam)/n^2+(-1/2\gam+9/4-1/2\ln(n))/n^3\\
&+3/(4n^4)+(1/20\gam+1/48+1/20\ln(n))/n^5+\BO(1/n^6).\\
\end{align*}
\section{The case $m=3,4$}\label{S8}
We have, for $m=3$, by (\ref{E1}),
\bals
Li_3(z)&=-1/2w^2(\ln(w)-3/2)+\zeta(3)-1/6\pi^2w+1/12w^3-1/288w^4+1/86400w^6-1/10160640w^8+\BO(w^9),\\
S_0(z)&=\zeta(3)/(1/2w^2(\ln(w)-3/2)+1/6\pi^2w-1/12w^3+1/288w^4-1/86400w^6+1/10160640w^8+\BO(w^9)),\\
&\mbox{we have the expansions}\\
w&=\eps+1/2\eps^2+1/3\eps^3+1/4\eps^4+\BO(\eps^5),\\
\ln(w)&=-L+1/2\eps+5/24\eps^2+1/8\eps^3+\BO(\eps^4).\\
&\mbox{Hence}\\
S(\eps,L)&=6\zeta(3)/(\pi^2\eps)+1/8\zeta(3)(-24\pi^6+144\pi^4L+216\pi^4)/\pi^8\\
&+1/8\zeta(3)(-4\pi^6-48\pi^4+432\pi^2L^2+1296\pi^2L+972\pi^2)/\pi^8\eps\\
&+1/8\zeta(3)(-2\pi^6-19\pi^4+360\pi^2L+216\pi^2L^2+54\pi^2+1296L^3+5832L^2+8748L+4374)/\pi^8\eps^2+\BO(\eps^3).\\
&\mbox{This leads successively to }\\
D_{-1,0}&=6\zeta(3)/\pi^2,\\
D_{0,1}&=18\zeta(3)/\pi^4,\\
D_{1,1}&=162\zeta(3)/\pi^6,\\
D_{1,2}&=54\zeta(3)/\pi^6,\\
T_{-1,0,n}&=6\zeta(3)/\pi^2,T_{-1,k,n}=0,k>0,\\
T_{0,1,n}&=18\zeta(3)/(\pi^4 n),T_{0,k,n}=0,k>1,\\
&\mbox{The general form of $G_{i,j}$ is given in  \cite{FlSe09}, Equ. (27). The detailed computation goes as follows},\\
G_{1,1,n}&=L_{1,n}-L_{1,n-1}=-1/n^2-1/n^3-1/n^4-1/n^5+\BO(1/n^6),\\
T_{1,1,n}&=D_{1,1}G_{1,1,n}=-162\zeta(3)/(\pi^6n^2)-162\zeta(3)/(\pi^6 n^3)+\BO(1/n^4),\\
G_{1,2,n}&=L_{2,n}-L_{2,n-1}=(2-2\ln(n)-2\gam)/n^2+(5-2\ln(n)-2\gam)/n^3\\
&+(43/6-2\ln(n)-2\gam)/n^4+(55/6-2\ln(n)-2\gam)/n^5+\BO(1/n^6),\\
T_{1,2,n}&=D_{1,2}G_{1,2,n}=54(2-2\ln(n)-2\gam)\zeta(3)/(\pi^6n^2)+54(5-2\ln(n)-2\gam)\zeta(3)/(\pi^6n^3)+\BO(1/n^4),\\
D_{2,3}& = 162\zeta(3)/\pi^8,\\
D_{2,2}&= 27\zeta(3)(27+\pi^2)/\pi^8,\\
D_{2,1}&= 9/2\zeta(3)(10\pi^2+243)/\pi^8,\\
G_{2,3,n}&=L_{3,n}-2L_{3,n-1}+L_{3,n-2}=(6-18\gam-\pi^2-18\ln(n)+6\gam^2+6\ln(n)^2+12\ln(n)\gam)/n^3+\BO(1/n^4),\\
G_{2,2,n}&=L_{2,n}-2L_{2,n-1}+L_{2,n-2}=(-6+4\ln(n)+4\gam)/n^3+\BO(1/n^4),\\
G_{2,1,n}&=L_{1,n}-2L_{1,n-1}+L_{1,n-2}=2/n^3+\BO(1/n^4),\\
T_{2,3,n}&=D_{2,3}G_{2,3,n}=162(6-18\gam-\pi^2-18\ln(n)+6\gam^2+6\ln(n)^2+12\ln(n)\gam)\zeta(3)/(n^3\pi^8)+\BO(1/n^4),\\
T_{2,2,n}&=D_{2,2,n}G_{2,2,n}=27(-6+4\ln(n)+4\gam)\zeta(3)(27+\pi^2)/(n^3\pi^8)+\BO(1/n^4),\\
T_{2,1,n}&=D_{2,1}G_{2,1,n}=9\zeta(3)(10\pi^2+243)/(n^3 \pi^8)+\BO(1/n^4).\\
&\mbox{Finally}\\
T_n&=T_{-1,0,n}+T_{0,1,n}+\sum_i\sum_j T_{i,j,n}.
\end{align*}
This leads to the following theorem:
\bth
\bals
&\mbox{Let}\\
 S_0(z)&:=\frac{Li_3(1)}{Li_3(1)-Li_3(z)},\mbox{ then}\\
S_n&:=[z^n]S_0(z) =6\zeta(3)/\pi^2+18\zeta(3)/(\pi^4n)+3\zeta(3)(-18\pi^2-36\pi^2\ln(n)-36\pi^2\gam)/(\pi^8n^2)\\
&+3\zeta(3)(-42\pi^2-405+324\gam^2+324\ln(n)^2+648\ln(n)\gam)/(\pi^8n^3)+\BO(1/n^4).
\end{align*}
\ethGL
\bals
&\mbox{This gives }\\
C_{n,0}&=6\zeta(3)/\pi^2,\\
C_{n,1}&=6\zeta(3)/\pi^2+18\zeta(3)/(\pi^4n),\\
C_{n,2}&=6\zeta(3)/\pi^2+18\zeta(3)/(\pi^4n)+3\zeta(3)(-18\pi^2-36\pi^2\ln(n)-36\pi^2\gam)/(\pi^8n^2),\\
C_{n,3}&=6\zeta(3)/\pi^2+18\zeta(3)/(\pi^4n)+3\zeta(3)(-18\pi^2-36\pi^2\ln(n)-36\pi^2\gam)/(\pi^8n^2)\\
&+3\zeta(3)(-42\pi^2-405+324\gam^2+324\ln(n)^2+648\ln(n)\gam)/(\pi^8n^3).
\end{align*}

\begin{figure}[htbp]
	\centering
		\includegraphics[width=0.8\textwidth,angle=0]{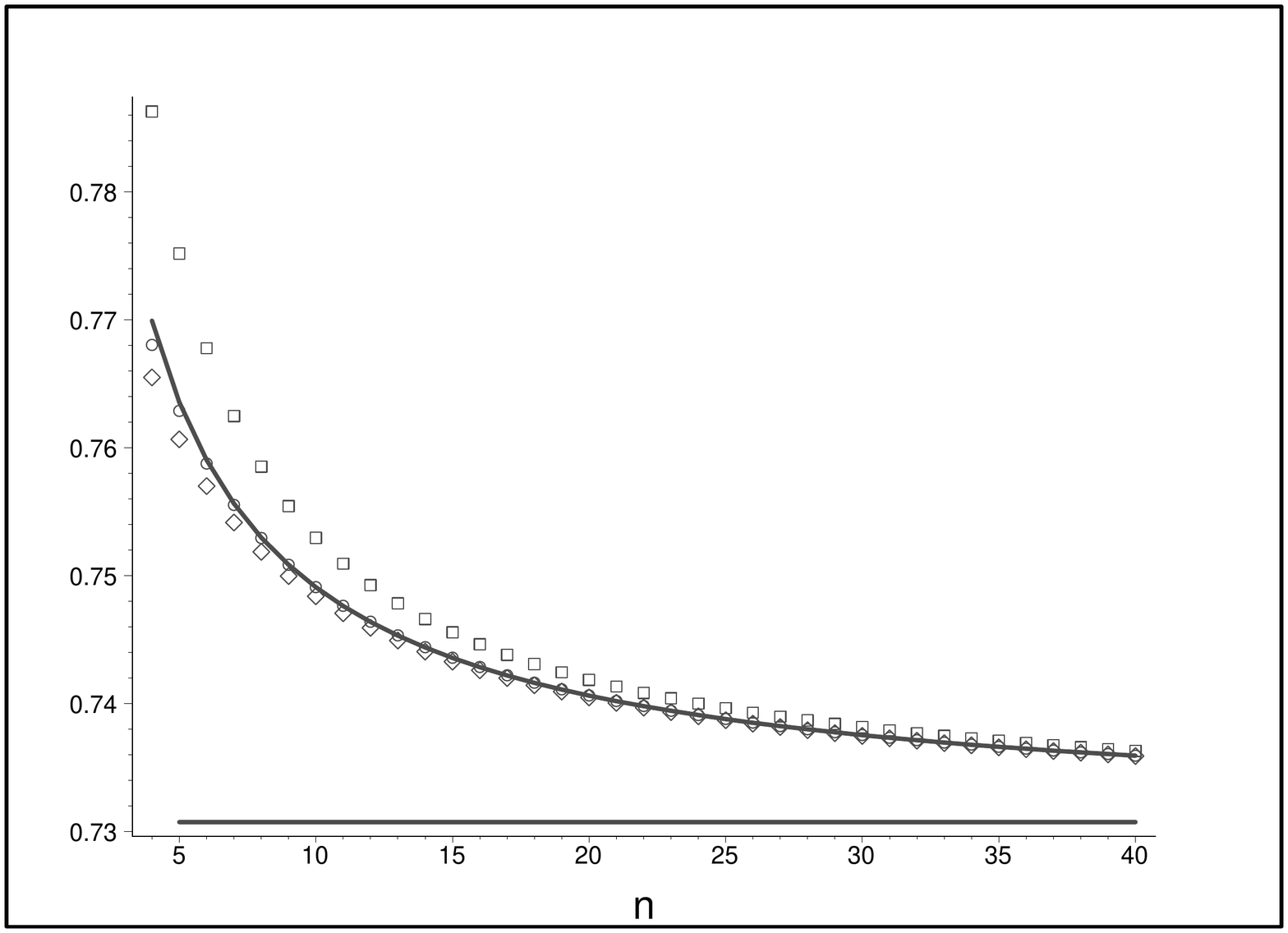}
	\caption{$m=3$, $S_n$(line),$C_{n,0}$(line),$C_{n,1}$(box),$C_{n,2}$(diamond),$C_{n,3}$(circle)}
	\label{F2}
\end{figure}

\newpage
To check the quality of our asymptotics, we display, in Figure \ref{F2}, $S_n,C_{n,0},C_{n,1},C_{n,2},C_{n,3}$.

The convergence of $S_n$ to $C_{n,0}$ is rather slow: we have $C_{n,0}=0.7307629692\ldots$ and $S_{100}=0.7329\ldots$.

\begin{figure}[htbp]
	\centering
		\includegraphics[width=0.8\textwidth,angle=0]{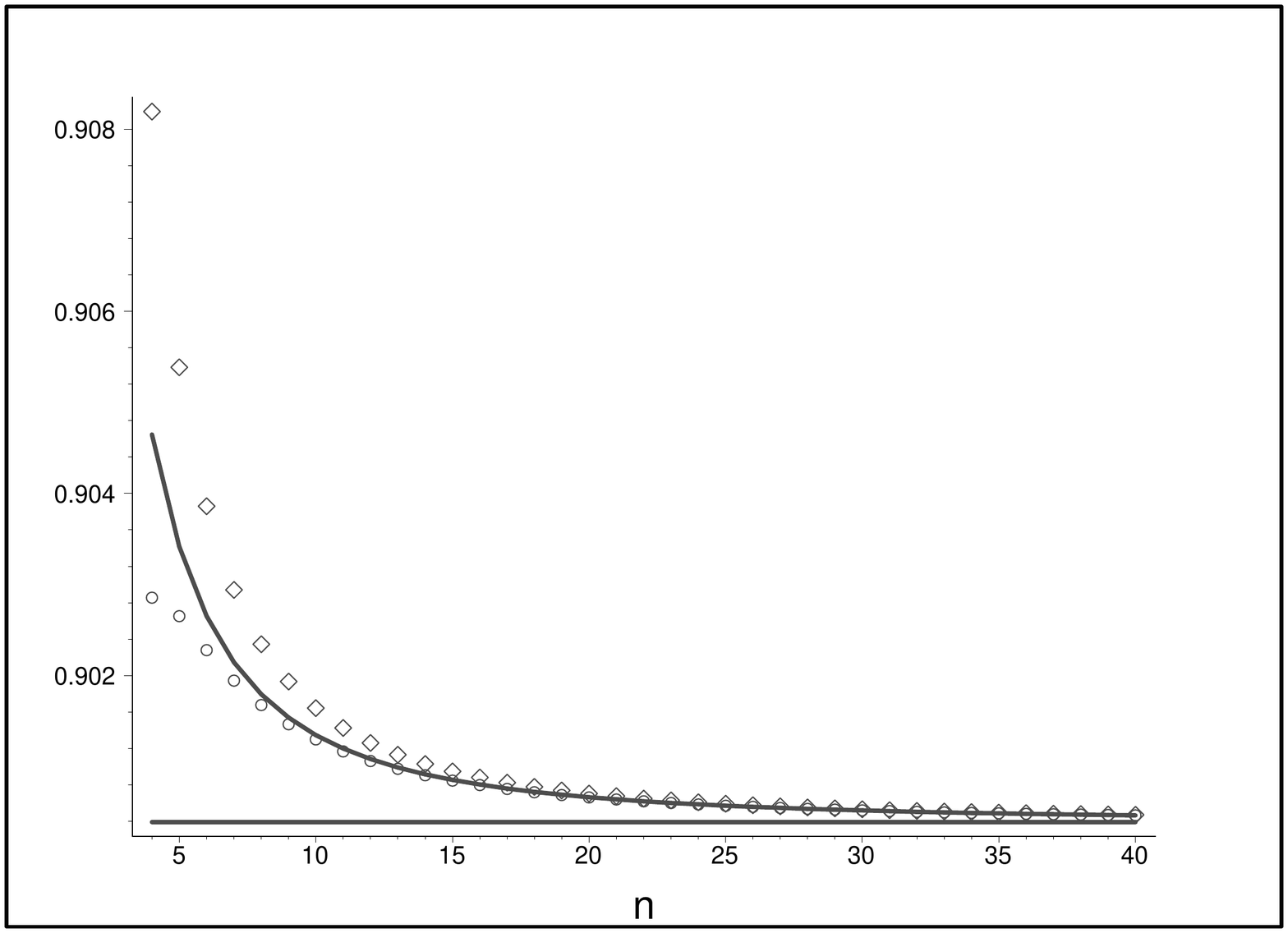}
	\caption{$m=4$, $S_n$(line),$C_{n,0}$(line),$C_{n,2}$(diamond),$C_{n,3}$(circle)  }
	\label{F3}
\end{figure}

The case $m=4$ is mechanically treated like the case $m=3$. We obtain
\bals
C_{n,0}&=1/90\pi^4/\zeta(3),\\
T_{0,1,n}&=0,\\
C_{n,2}&=1/90\pi^4/\zeta(3)+1/540\pi^4/(\zeta(3)^2n^2),\\
C_{n,3}&=1/90\pi^4/\zeta(3)+1/540\pi^4(\zeta(3)^2n^2)-1/1620\pi^6/(\zeta(3)^3 n^3),
\end{align*}
we display, in Figure \ref{F3}, $S_n,C_{n,0},C_{n,2},C_{n,3}$.
The convergence of $S_n$ to $C_{n,0}$  is faster.
\vspace{1em}

Remarks

\ben
\item in order to expand $S_0(1-\eps)$, we first expand w.r.t $\eps$ as $\eps=o(1/L^k),k>0$
\item $G_{k,j}$ starts with a $1/n^{k+1}$ term, which allows an easy expansion
\item  more and more terms  in the expansion of $S(\eps,L)$ are needed when $m$ increases: the first terms 
don't contain any $L$ terms. For instance, for $m=6$, only the $\eps^3$ contains a linear $L$ contribution,
and the asymptotics of  $S_n$ starts as $1/945\pi^6/\zeta(5)+1/18900  \pi^6/(\zeta(5)^2n^4)$.  More terms can be mechanically computed.
\een
\section{Acknowledgements}\label{S9}
We would like to thank H.Prodinger for helping in computing an Euler sum  and W.Wang for providing a useful reference.
\bibliographystyle{plain}

\begin{thebibliography}{1}

\bibitem{AP10}
T.M. Apostol.
\newblock {\em Polylogarithm}.
\newblock in NIST Handbook of Mathematical Functions. Cambridge University
  Press, 2010.

\bibitem{DO91}
P.J. de~Doelder.
\newblock On some series containing $\psi(x)-\psi(y)$ and $(\psi(x)-\psi(y))^2$
  for certain values of $x$ and $y$.
\newblock {\em Journal of Computational and Applied Mathematics},
  37(1-3):125--141, 1991.

\bibitem{FlSa98}
P.~Flajolet and B.~Salvy.
\newblock Euler sums and contour integral representations.
\newblock {\em Journal of Experimental Mathematics}, 7(1):15--35, 1998.

\bibitem{FlSe09}
P.~Flajolet and R.~Sedgewick.
\newblock {\em Analytic Combinatorics}.
\newblock Cambridge University Press, Cambridge, 2009.

\bibitem{GOWA17}
R.~G\'{o}mez-Aiza and M.D. Ward.
\newblock 2017.
\newblock private communication.

\bibitem{LE81}
L.~Lewin.
\newblock {\em Polylogarithms and associated functions}.
\newblock North-Holland, 1981.

\bibitem{XU17}
C.~Xu.
\newblock Evaluations of {E}uler type sums of weight $\leq 5$.
\newblock Technical report, School of Mathematical Sciences, Xiamen University,
  2017.
\newblock arXiv preprint arXiv:1704.03515.

\end{thebibliography}

\end{document}